\documentclass[11pt]{article}

\usepackage{amsfonts,latexsym,amstext}
\usepackage{amsmath,amscd,amssymb}
\usepackage[active]{srcltx}

\usepackage{xcolor}

\newtheorem{theorem}{Theorem}[section]
\newtheorem{lemma}[theorem]{Lemma}
\newtheorem{proposition}[theorem]{Proposition}
\newtheorem{definition}{Definition}[section]
\newtheorem{example}{Example}[theorem]
\newtheorem{hypothesis}[theorem]{Hypothesis}

\newtheorem{remark}[theorem]{Remark}
\newtheorem{corollary}[theorem]{Corollary}
\newenvironment{proof}[1][Proof]{\medskip \noindent \textbf{#1.}
 }{\ \rule{0.5em}{0.5em}}
\makeatletter \@addtoreset{equation}{section}

\makeatother

\def\qed{{\hfill\hbox{\enspace${ \square}$}} \smallskip}
\def\sqr#1#2{{\vcenter{\vbox{\hrule height .#2pt \hbox{\vrule
 width .#2pt height#1pt \kern#1pt \vrule
width .#2pt} \hrule height .#2pt}}}}
\def\square{\mathchoice\sqr54\sqr54\sqr{4.1}3\sqr{3.5}3}

\def\ds{\begin{displaystyle}}
\def\eds{\end{displaystyle}}

\def\<{\langle }
\def\>{\rangle }

\def\R{\mathbb R}

\def\E{\mathbb E}
\def\P{\mathbb P}

\def\calf{{\cal F}}
\def\calg{{\cal G}}

\def\calp{{\cal P}}

\title{Ergodic BSDEs under weak dissipative assumptions}

\author{Arnaud Debussche\\ENS de Cachan, Antenne de Bretagne\\Campus de Ker Lann, Av. R. Schuman, 35170 Bruz,
France\\arnaud.debussche@bretagne.ens-cachan.fr\\ \\
Ying Hu\\IRMAR, Universit\'e Rennes 1\\Campus de Beaulieu, 35042
Rennes Cedex, France\\ying.hu@univ-rennes1.fr\\ \\
Gianmario
Tessitore\\Dipartimento di Matematica e Applicazioni, Universit\`a
di Milano-Bicocca\\Via Cozzi 53, 20125 Milan,
Italy\\gianmario.tessitore@unimib.it}

\begin{document}

\maketitle

\begin{abstract}
In this paper we study ergodic backward stochastic differential equations (EBSDEs) dropping the strong dissipativity
assumption needed in \cite{FHT-Erg}. In other words we do not  need to require the uniform exponential decay of the difference of two solutions of the underlying forward equation, which, on the contrary, is assumed to be non degenerate.

We show existence of solutions
by use of coupling estimates for a non-degenerate forward stochastic differential equations
with bounded measurable non-linearity. Moreover we prove uniqueness of ``Markovian'' solutions
exploiting the recurrence of the same class of forward equations.

Applications are then given to the optimal ergodic control of stochastic partial differential equations
 and to the associated ergodic Hamilton-Jacobi-Bellman equations.
\end{abstract}

\section{Introduction}
Since the beginning of the 90's several papers have described the link between
 backward stochastic differential equations (BSDEs), Hamilton-Jacobi-Bellman equations
 and stochastic optimal control (see, for instance, \cite{kapequ} and \cite{peng93}).
The successive literature on BSDEs covered several different situations, among
 them  infinite horizon control problems,  both in finite and infinite dimensions
(see \cite{bh}, \cite{royer} and \cite{HuTe}).

In \cite{FHT-Erg} the BSDE approach was extended to the case of ergodic
 control problems, that is, of control problems in which the
cost functional only evaluates  the long time behavior of the stochastic system.
In that paper the authors introduced the following class of BSDEs with infinite horizon, called Ergodic BSDEs (EBSDEs):
\begin{equation}\label{ebsde}
Y_t^x=Y_T^x+\int_t^T [\psi(X^x_\sigma,Z^x_\sigma)-\lambda]d\sigma -\int_t^T Z_\sigma^xdW_\sigma,
\quad \mathbb{P}-\hbox{a.s.},\quad \forall\, 0\le t\le T<\infty,
\end{equation}
where $(W_t)_{t\ge 0}$ is a cylindrical Wiener process in a Hilbert space $\Xi$, $X^x$ is the solution
of the following forward SDE with values in a Hilbert (in \cite{FHT-Erg} also Banach) space $H$,
\begin{equation}\label{forward-intro} dX_t^x=(AX_t^x+F(X_t^x))dt+GdW_t,\quad X_0=x,\end{equation}
and $\psi:H\times \Xi^*\rightarrow \mathbb{R}$ is a given function.

We underline that the unknowns in the above equation  is the triple $(Y,Z,\lambda)$,
 where $Y,Z$ are adapted processes taking values in $\mathbb{R}$ and $\Xi^*$,
respectively, and $\lambda$ is a real
number.

The case of the EBSDEs driven by a finite dimensional reflected
forward equation together with its applications to semilinear PDEs
with general Neumann boundary conditions was then treated in
\cite{richou}.

The main assumption in \cite{FHT-Erg} (and, with slight
modifications, in \cite{richou}) is the strict dissipativity of
$A+F$, i.e., the existence of $k>0$ such that for all $x$ and $x'$ in the
domain of $A$, the following holds
$$(A(x-x')+F(x)-F(x'),x-x')\le -k |x-x'|_H^2.$$
Such a requirement ensures the uniform exponential decay of the difference
between the trajectories of two solutions of equation (\ref{forward-intro})
that plays a crucial role in the  arguments in \cite{FHT-Erg}.

$ $

The aim of the present paper is to show that, when  $G$ is invertible, we can drop the dissipativity assumption on $A+F$ and
study (\ref{ebsde})  when $A$ is
dissipative but $F$ is only assumed to be  bounded and Lipschitz with
no restrictions on its Lipschitz constant. See Example \ref{ex-heat-eq} to compare,
in the concrete case of an ergodic problem for a stochastic heat equation, the assumptions needed in the present paper
 and the ones needed in \cite{FHT-Erg} .

Our main tool is a  coupling estimate for  a perturbed version of the
forward stochastic differential equation (\ref{forward-intro}).
Coupling estimates have been recently developed  for many different classes of stochastic partial differential equations
 and exploited to deduce regularity properties of the corresponding Markov semigroup,
see e.g. \cite{EMS}, \cite{kuksin-shirikyan}, \cite{mattingly},
\cite{shirikyan} and  \cite{odasso}. In the present paper, in comparison with the previously mentioned literature,
we are  only dealing with bounded and everywhere defined non-linearities
but we have to consider measurable non-linearities and to  prove
that the  estimate depends on them only through their sup (see Theorem \ref{coup} and Appendix \ref{sec-app-coupl}).

The coupling estimate is  used here to get a uniform bound for ${Y}_t^{x,\alpha}-{Y}^{0,\alpha}_0$ where
${Y}^{x,\alpha}$ is the solution of the following strictly monotonic BSDE (see \cite{royer})
$${Y}^{x,\alpha}_t={Y}^{x,\alpha}_T +\int_t^T(\psi(X^{x}_\sigma,
Z^{x,\alpha}_\sigma)-\alpha Y^{x,\alpha}_\sigma)d\sigma-\int_t^T
Z^{x,\alpha}_\sigma dW_\sigma, \quad 0\le t\le T <\infty. $$
Then the Bismut-Elworthy formula for BSDEs (see \cite{FuTeBE}) yields the uniform bound for
$\nabla_x {Y}^{x,\alpha}_0$ that allows to pass to the limit as $\alpha \searrow 0$  in the above equation. Note that the non degeneracy of the noise, i.e., the invertibility of $G$, is used in an 
essential way at this step. It is also used to prove the coupling estimate but a more sophisticated
coupling argument, which would not need this assumption,  could be used. 

We also notice that we construct a ``Markovian'' solution of the EBSDE  in the sense that $Y_t$ and $Z_t$
are deterministic functions of $X^x_t$.
We prove,  by use of the recurrence of the perturbed forward stochastic differential equation, that  such a ``Markovian''
solution of the EBSDE is unique. The recurrence property is studied in \cite{DP2} for a forward SDE similar to ours;
the difference is that here we need to consider drifts that are only bounded and measurable.
 (see Theorem 2.6 and Appendix \ref{sec-app-rec}).
The uniqueness argument is inspired by the corresponding one in \cite{GoMa}.

$ $

Once existence and uniqueness of a Markovian solution of the EBSDE is proved we can proceed as in \cite{FHT-Erg}
to deal with an optimal control problem with state equation
\begin{equation}\label{se-intro}dX_t^{x,u}=(AX_t^{x,u}+F(X_t^{x,u})+ GR(u_t))dt+GdW_t,\quad X^{x,u}_0=x,\end{equation}
and ergodic cost functional
\begin{equation}\label{-ergodic-cost-intro}
  J(x,u)=\limsup_{T\rightarrow\infty}\frac{1}{T} \mathbb
E^{u,T}\int_0^T L(X_s^x,u_s)ds.
\end{equation}
Then we deduce that the ergodic Hamilton Jacobi Bellman equation
 \begin{equation}
\mathcal{L}v(x)
+\psi\left( x,\nabla v(x)  G\right)  = \lambda, \quad x\in E,  \label{hjb*}%
\end{equation}
has a unique mild solution; moreover the  ergodic problem admits a unique optimal control that satisfies
 an optimal feedback law
given in terms of the gradient on the solution to the  HJB equation (\ref{hjb*}); finally the optimal  cost is  $\lambda$.

$ $

In the finite dimensional case there are several papers devoted to
the study, by analytic techniques,  of stochastic optimal ergodic control problems
 and of the corresponding HJB equations (see for instance
 \cite{ArLi} and  \cite{BeFr}). On the contrary, to the best of our knowledge,
 there are very few works devoted to the infinite dimensional case.
As far as we know (mild) solutions of an equation like (\ref{hjb*})
was studied, in the infinite dimensional case, only in
\cite{GoMa} (besides the already discussed results included in \cite{FHT-Erg}).

In \cite{GoMa} authors  prove, by a fixed point argument, existence and
uniqueness of the  solution of the mild stationary HJB equation for
discounted infinite horizon costs. Then they pass to the limit, as
the discount goes to zero. They work under the same
non-degeneracy assumption that we use here and assume that $A$ is the generator
 of a contraction semigroup and $ F$ is dissipative; they also have
 a limitation on the Lipschitz constant (with respect to
the gradient variable)  of the Hamiltonian function $\psi$ (see
\cite{FuTe-ell} for similar conditions in the case of a  strictly monotonic stationary HJB equation ).
On the contrary unbounded non-linearities $F$ can be considered in \cite{GoMa}.

The present paper is organized as follows. First we establish general notation.
In section \ref{sec-forward}, we introduce the forward equation and state the coupling estimates and recurrence property
for the perturbed forward equation.
The ergodic BSDE is studied  in Section \ref{sec-ergodic}. Sections \ref{sec-HJB} and \ref{sec-contr} shortly recall
how the previous results can be applied to the ergodic Hamilton-Jacobi-Bellman equation and
 to the  ergodic optimal control problem. We include the proofs
for the coupling estimates and the recurrence property for the perturbed forward equation in the Appendix (see section \ref{appendix}).

\subsection{General Notation}

We introduce some notations; let $E,F$ be real separable Hilbert spaces. The
norms and the scalar product will be  denoted $|\,\cdot\,|$,
$\langle\,\cdot\,,\,\cdot\,\rangle$, with subscripts if needed.
$L(E,F)$ is the space of linear bounded operators $E\to F$, with
the operator norm.
The domain of a linear (unbounded) operator $A$ is denoted $D(A)$.

Given $\phi\in B_b(E)$, the space of bounded and measurable functions
$ \phi: E\rightarrow \mathbb{R}$, we denote
$\Vert\phi\Vert_0=\sup_{x\in E}|\phi(x)|$. If, in addition,
$\phi$ is also Lipschitz continuous then
$\Vert\phi\Vert_{\hbox{lip}}=\Vert\phi\Vert_0+
\sup_{x,x'\in E,\,x\ne x'}|\phi(x)-\phi(x')||x-x'|^{-1}$.

We say that a function $F:E\to F$ belongs to
the class $\calg^1(E,F)$ if it is continuous, has a Gateaux
differential $\nabla F(x)\in L(E,F)$ at any point $x  \in E$, and
for every $k\in E$ the mapping $x\to \nabla F(x) k$  is continuous
from $E$ to $F$ (i.e. $x\to \nabla F(x) $  is continuous from $E$
to $L(E,F)$ if the latter space is endowed with the strong operator
topology). In connection with stochastic equations,
the space $\calg^1$ has been introduced in \cite{FuTe1},
to which we refer the reader for further properties.

 Given a  probability space $\left(
\Omega,\mathcal{F},\mathbb{P}\right) $ with a filtration
$(\calf_t)_{t\ge 0}$ we consider the following classes of
stochastic processes with values in a real separable Banach space
$K$.

\begin{enumerate}
\item
$L^p_{\mathcal{P}}(\Omega,C([0,T],K))$, $p\in [1,\infty)$,
$T>0$, is the space
of predictable processes $Y$ with continuous paths
on $[0,T]$
such that
$$
|Y|_{L^p_{\mathcal{P}}(\Omega,C([0,T],K))}^p
= \E\, \sup_{t\in [0,T]}|Y_t|_K^p<\infty.
$$
\item
$L^p_{\mathcal{P}}(\Omega,L^2([0,T];K))$, $p\in [1,\infty)$,
$T>0$, is the space
of predictable processes $Y$ on $[0,T]$ such that
$$
|Y|^p_{L^p_{\mathcal{P}}(\Omega,L^2([0,T];K))}=
\E\,\left( \int_{0}^{T}|Y_t|_K^2\,dt\right)^{p/2}<\infty.
$$

\item
$L_{\cal P, {\rm loc}}^2(\Omega;L^2(0,\infty;K))$
is the space
of predictable processes $Y$ on $[0,\infty)$ that belong
to the space $L^2_{\mathcal{P}}(\Omega,L^2([0,T];K))$
for every $T>0$.
\end{enumerate}

\section{The forward SDE}\label{sec-forward}

\subsection{General assumptions}

This section is devoted to the following mild It\^o stochastic
differential equation for an unknown process $X_\tau$,
$\tau\in \mathbb{R}^+$, with values in a Hilbert space $H$:
\begin{equation}\label{equazioneforwardmild}
\hat{X}_\tau = e^{(\tau-t)A}x+\int_t^\tau e^{(\tau-\sigma)A}
\Upsilon(\hat{X}_\sigma)\; d\sigma
+\int_t^\tau e^{(\tau-\sigma)A}
G\, dW_\sigma,
\qquad \forall \tau\geq 0,\quad \mathbb{P}-\hbox{a.s.}.
\end{equation}
We assume the following:
\begin{hypothesis}\label{ipotesiuno} $ $

\begin{description}
  \item[(i)] $A$ is an unbounded operator $A:D(A)\subset H
\rightarrow H$, with $D(A)$ dense in $H$.
We assume that $A$ is  dissipative
and generates a stable $C_0$-semigroup $\{e^{tA}\}_{t\geq
0}$. By this we mean
 that there exist constants $k>0$ and $M>0$ such that
$$\langle A x,x \rangle \leq -k |x|^2\quad \forall x\in D(A);
\qquad |e^{\tau A}|\leq M e^{-k \tau}.$$

\item[(ii)] For all $s>0$, $e^{sA}$ is a Hilbert-Schmidt operator. Moreover
$ |e^{sA}|_{L_2(H,H)}\leq L\; s^{-\gamma}$
for  suitable constants $L>0$ and $\gamma\in [0,1/2)$.

\item[(iii)] $\Upsilon$ is a bounded measurable map $H\to H$,
\item[(iv)] $G$ is a bounded linear operator in $L(\Xi,H)$.
 Moreover we assume that $G$ is invertible and
 we denote by $G^{-1}$ its bounded inverse.

\item[(v)] $(\Omega,\mathcal{F},\mathbb{P})$ is a complete probability space, $(\mathcal{F}_t)_{t\geq 0}$ is a filtration in it satisfying the usual conditions and $(W_t)_{t\geq 0}$ is an $\mathcal{F}$-cylindrical Wiener process with values in a separable Hilbert space $\Xi$.

  \end{description}
\end{hypothesis}

\begin{remark}\label{A_-m-dissip}{\em We notice that if the operator
$A$ with dense domain is m-dissipative that is $\langle A x,x
\rangle \leq -k |x|_H^2\quad \forall x\in D(A)$ and $A-k_1 I$ is
surjective for a suitable $k_1>0$ then by the Lumer-Phillips theorem it follows
immediately  that $A$ generates stable $C_0$-semigroup of
contractions (that is $M=1$).}
\end{remark}


The following result is well known in its first part (see, for instance, \cite{DP1}) and a straight-forward consequence
of the Girsanov transform in the second.

\begin{proposition}\label{esisteunicox} Fix $t\geq 0$ and $x\in H$ and assume that $\Upsilon$ is Lipschitz.
Under the assumptions of Hypothesis
 \ref{ipotesiuno} there
exists a unique  adapted process
 $\hat{X}$ verifying (\ref{equazioneforwardmild}). Moreover, for every $p\in [2,\infty)$ and every $T>t$,
 $\hat{X} \in L_{\calp}^p(\Omega;C([0,T];H))$ and
\begin{equation}\label{stimadeimomentidix}
  \E\sup_{\tau\in[t,T]}|X_\tau |^p\leq C(1+|x|)^p,
\end{equation}
for some constant $C$ depending only on $p,\gamma,M$ and
$\sup_{x\in H}|\Upsilon(x)|$ but independent of $T>t$.

$ $

If $\Upsilon$ is only bounded and  measurable,  then the solution to equation (\ref{equazioneforwardmild})
still exists but in \emph{weak} sense. By this we mean, see again \cite{DP1}, that there exists a new $\mathcal{F}$-Wiener process $(\hat{W}_t)_{t\geq 0}$  with respect to a new probability
$\hat{\mathbb{P}}$ (absolutely continuous with respect to $\mathbb{P}$),
 and an $\mathcal{F}$-adapted process  $\hat{X}$ with continuous trajectories
for which (\ref{equazioneforwardmild}) holds
with $W$ replaced by $\hat W$.
Moreover (\ref{stimadeimomentidix}) still holds (with respect to the new probability).
Finally such a weak solution is unique in law.
\end{proposition}

In the following we will denote  the solution of equation (\ref{equazioneforwardmild}) by $\hat{X}^{t,x}$   and by $\hat{X}^{x}$ when we
choose the initial time $t=0$. We remark that equation (\ref{equazioneforwardmild})
 is the mild version of the Cauchy problem:
\begin{equation}
\left\{
\begin{array}[c]{l} d\hat X^{t,x}_s  =A\hat X^{t,x}_s
 d s+\Upsilon(\hat X^{t,x}_s) dt +G dW_s ,\text{ \ \ \
} s\geq t, \\
\hat X^{t,x}_t  =x.
\end{array}
\right.  \label{sde-coup}
\end{equation}
The following result is proved in section 6.1.
\begin{theorem}[Basic coupling estimate]\label{coup}
Assume that $\Upsilon: H\rightarrow H$ is Lipschitz and let $\hat X^x$ be the (strong)
solution of equation (\ref{equazioneforwardmild})
then there exist $\hat c>0$  and $\hat \eta>0$ such that for all
$\phi \in B_b(H)$
\begin{equation}  \label{coupling-estimate}
\left|\mathcal{P}_t[\phi](x) - \mathcal{P}_t[\phi](x')\right| \leq
\hat c(1+|x|^2+|x'|^2) e^{-\hat \eta t} |\phi|_0,
\end{equation}
where $\mathcal{P}_t[\phi](x)=\mathbb{E} \phi (\hat X^x_t )$ is the
Kolmogorov semigroup associated to equation (\ref{equazioneforwardmild}).

 We stress  the fact that  $\hat c$ and $\hat \eta$ depend on
$\Upsilon$ only through $\sup_{x\in H}|\Upsilon(x)|$.
 \end{theorem}

 \begin{corollary}\label{cor-coupl} Relation (\ref{coupling-estimate})
can be extended to the case in which $\Upsilon$ is only bounded and measurable
 and there exists a uniformly bounded sequence of Lipschitz functions
$\{\Upsilon_n\}_{n\ge 1}$ (i.e. $\forall n, \Upsilon_n$ is Lipschitz
and $\sup_n\sup_x|\Upsilon_n(x)|<\infty$) such that
$$\lim_n \Upsilon_n(x)=\Upsilon(x),\quad \forall x\in H.$$
Clearly in this case in the definition of $\mathcal{P}_t[\phi]$ the mean
value is taken with respect to the new probability $\hat P$.

\end{corollary}
\proof
It is enough to show that if $\mathcal{P}^{n}$ is the semigroup
corresponding to equation (\ref{equazioneforwardmild}) with $\Upsilon$ replaced
by $\Upsilon_n$, then $\forall x\in H$ and $\forall t\geq 0$,
$$ \mathcal{P}^{n}_t[ \phi] (x)\rightarrow \mathcal{P}_t[ \phi] (x).$$

We set
$$U^x_{\tau}=e^{\tau A}x
+\int_0^\tau e^{(\tau-\sigma)A}
G\; dW_\sigma.$$

By Girsanov's formula $$\mathcal{P}^{n}_t[ \phi]
(x)=\mathbb{E}(\rho^{n,x}_t \phi(U^x_t)),\qquad \mathcal{P}_t[ \phi]
(x)=\mathbb{E}(\rho^{x}_t \phi(U^x_t))$$  where
$$ \rho^{n,x}_t=\exp\left( -\int_0^t \<G^{-1}\Upsilon_n(U^x_t), dW_s\>_{\Xi}-
\frac{1}{2} \int_0^t |G^{-1}\Upsilon_n(U^x_t)|^2_{\Xi} ds\right),$$
and
$$ \rho^{x}_t=\exp\left( -\int_0^t \<G^{-1}\Upsilon(U^x_t), dW_s\>_{\Xi}-
\frac{1}{2} \int_0^t |G^{-1}\Upsilon(U^x_t)|^2_{\Xi} ds\right).$$

We have
\begin{eqnarray*}
& &\mathbb{E}[(\rho^{n,x}_t)^2]\\
&=&\mathbb{E}\left[\exp\left( -2\int_0^t\!\! \<G^{-1}\Upsilon_n(U^x_t), dW_s\>_{\Xi}-
2 \int_0^t\!\! |G^{-1}\Upsilon_n(U^x_t)|^2_{\Xi} ds\right)\exp\left( \int_0^t\!\!
 |G^{-1}\Upsilon_n(U^x_t)|^2_{\Xi} ds\right)\right]\\
&\le&\exp\left( t |G^{-1}|^2\sup_n\sup_x|\Upsilon_n(x)|^2 \right)<\infty,
\end{eqnarray*}
from which we deduce that $\{\rho^{n,x}_t\}_n$ is uniformly integrable in $L^{1}(\Omega)$.
Moreover, it is easy to see that $\lim_n \rho^{n,x}_t=\rho^{x}_t$ in probability, and the claim follows. \qed


$ $

The equation (\ref{equazioneforwardmild}) also enjoys a recurrence property that will be useful in the following, it is proved in section 6.2.

\begin{theorem}\label{recurrence}
Assume that $\Upsilon: H\rightarrow H$  can be approximated (in the sense of pointwise convergence)
by a uniformly bounded sequence of Lipschitz functions $\{\Upsilon_n\}_{n\ge 1}$.
Then the solution of equation (\ref{equazioneforwardmild}) is recurrent in the sense that
 for all $\Gamma\in H$, $\Gamma $ open:
$$\lim_{T\to\infty}\hat{\mathbb{P}}\{\exists t\in[0,T] : \hat X^x_t\in \Gamma\}=1.$$
In particular, setting
 $\tau^x =\inf\{ t : |\hat X^x_t|<
\epsilon\}$,  then
$\forall \epsilon>0$, $\lim_{T\rightarrow\infty}\hat{ \mathbb P}\{ \tau^x<T \}=1.$
\end{theorem}

\section{The Ergodic BSDE}\label{sec-ergodic}

We fix now a bounded function $F:H\rightarrow H$ and denote by
 ${X}^{t,x}$ (and by ${X}^{x}$ when we
choose the initial time $t=0$) the solution of equation (\ref{equazioneforwardmild})
 with $\Upsilon=F$.

\noindent This section is devoted to the following type of BSDEs with
infinite horizon
\begin{equation}\label{EBSDE}
Y^x_t=Y^x_T +\int_t^T\left[\psi(X^x_\sigma,Z^x_\sigma)-
\lambda\right]d\sigma-\int_t^T Z^x_\sigma\, dW_\sigma, \quad 0\le
t\le T <\infty,
\end{equation}
where $\lambda$ is a real number and is part of the unknowns of
the problem; the equation is required to hold for every $t$ and
$T$ as indicated. On the function $\psi: H\times  \Xi^*
\rightarrow {\mathbb R}$ and $F$ we assume the following:

\begin{hypothesis}\label{hypothesisroyer} $\psi$ is a measurable map $H\times  \Xi^*
\rightarrow {\mathbb R}$. Moreover there exists $l>0$  such that
$$ |\psi(x,0)|\le l; \quad |\psi(x,z)
-\psi(x,z')|\le l|z-z'|, \qquad
 x\in H,\;
z,z'\in\Xi^*.
$$
\end{hypothesis}
\begin{hypothesis}\label{hypothesis_F}  $F$ is bounded,  Lipschitz and  G\^ateaux differentiable,
more precisely,  $F$ belongs to the class
$\calg^1(H, H). $
\end{hypothesis}
We start by considering an infinite horizon equation with strictly
monotonic drift, namely, for $\alpha>0$, the equation
\begin{equation}\label{bsderoyer}
{Y}^{x,\alpha}_t={Y}^{x,\alpha}_T +\int_t^T(\psi(X^{x}_\sigma,
Z^{x,\alpha}_\sigma)-\alpha Y^{x,\alpha}_\sigma)d\sigma-\int_t^T
Z^{x,\alpha}_\sigma dW_\sigma, \quad 0\le t\le T <\infty.
\end{equation}

The existence and uniqueness of solution to (\ref{bsderoyer})
under Hypothesis \ref{hypothesisroyer} was first studied by Briand
and Hu in \cite{bh} and then generalized by Royer in \cite{royer}.
The following lemma follows from Lemma 2.1 and Theorem 3.2 in \cite{HuTe}.

\begin{lemma}\label{lemmaroyer} Let us suppose that the Hypotheses
\ref{ipotesiuno},
\ref{hypothesisroyer} and 3.2  hold.
 Then for all $x\in H$ and $\alpha>0$,

 (i) there  exists a unique solution $(Y^{x,\alpha},Z^{x,\alpha})$
 to  BSDE (\ref{bsderoyer})
such that $Y^{x,\alpha}$ is a bounded continuous process,
$Z^{x,\alpha}$ belongs to $L_{\cal P, {\rm
loc}}^2(\Omega;L^2(0,\infty;\Xi^*))$, and $|Y^{x,\alpha}_t|\leq {l}/{\alpha}$, $\mathbb{P}$-a.s.
for all $t\geq 0$;

(ii)  if we define $v^{\alpha}(x)=Y^{x,\alpha}_0$ then, for all fixed $\alpha>0$,  $v^{\alpha}$ is Lipschitz  bounded and of class ${\cal G}^1$, moreover,
$$Y_t^{x,\alpha}=v^\alpha(X_t^x),\quad Z_t^{x,\alpha}=\nabla v^\alpha(X_t^x)G.$$
\end{lemma}

In order to construct the solution to (\ref{EBSDE}), we need some uniform in $\alpha$ estimate of $|v^\alpha(x)-v^\alpha(x')|$. This will be obtained by coupling estimates but first we have to prove  an approximation lemma:
\begin{lemma}\label{linear}
Let $\zeta,\; \zeta': H \rightarrow \Xi^*$ weakly$^*$ continuous
with polynomial growth. We define
$$\tilde{\Upsilon}(x)=\begin{cases} \dfrac{\psi(x,\zeta(x)) -\psi(x,\zeta'(x))}{|\zeta(x))
 -\zeta'(x)|^2}
  \left(\zeta(x)-\zeta'(x)\right)^*, &  \hbox{ if $\; \zeta(x)
 \neq \zeta'(x)$}, \\ \\
  0,  \hbox{ if $\; \zeta(x)
= \zeta'(x).$}
 \end{cases}  $$
There exists a uniformly bounded sequence of Lipschitz functions $(\tilde{\Upsilon}_n)_{n\ge 1}$
 (i.e., $\forall n$, $\tilde{\Upsilon}_n$ is Lipschitz and
$\sup_n\sup_x |\tilde{\Upsilon}_n(x)|<\infty$) such that
$$\lim_n \tilde{\Upsilon}_n(x)=\tilde{\Upsilon}(x),\quad \forall x\in H.$$

\end{lemma}
\proof Fixing an orthonormal basis $\{\xi_1,\xi_2,\cdots\}$ in $\Xi$, we define  the projection $\Pi_{p,\Xi^*}$: $\Xi^*\rightarrow
\Xi^*$ as follows:
$$\Pi_{p,\Xi^*}\zeta=\sum_{i=1}^n (\zeta\xi_i)<\xi_i,\cdot>.$$
Set
$$ \widetilde{\Upsilon}^{i}(x)=\dfrac{\psi(x,\zeta(x)) -\psi(x,\zeta'(x))}{|\zeta(x)
 -\zeta'(x)|^2+i^{-1}}
  \left(\zeta(x)-\zeta'(x)\right)^*, $$
  $$ \widetilde{\Upsilon} ^{i,p}(x)=\dfrac{\psi(x,\Pi_{p,\Xi^*}\zeta(x))
  -\psi(x,\Pi_{p,\Xi^*}\zeta'(x))}{|\Pi_{p,\Xi^*}\zeta(x))
 -\Pi_{p,\Xi^*}\zeta'(x)|^2+i^{-1}}
  \left(\Pi_{p,\Xi^*}\zeta(x)-\Pi_{p,\Xi^*}\zeta'(x)\right)^*. $$
  It is easy to verify that the functions $\widetilde{\Upsilon}
  ^{i,p}$ are continuous functions.  Moreover
  $|\widetilde{\Upsilon} ^{i,p}(x)| \leq l$ and
$ \lim_p\widetilde{\Upsilon} ^{i,p}(x)= \widetilde{\Upsilon}
^{i}(x)$ , $\lim_i\widetilde{\Upsilon} ^{i}(x)=
\widetilde{\Upsilon} (x)$, for all $x\in H$. Fixing $i,p$, it is quite classical
(based on finite dimensional projections and convolutions) to construct a uniformly bounded sequence of Lipschitz functions
$\{\widetilde{\Upsilon} ^{i,p,m}\}_m$, such that $\lim_m\widetilde{\Upsilon} ^{i,p,m}(x)=\widetilde{\Upsilon} ^{i,p}(x)$,
 see, e.g. Lemma 4.2 in \cite{FuTeBE}. Then the proof ends with a diagonal procedure.
 \qed

$ $

The following lemma
plays a crucial role. It gives the desired estimate of $v^\alpha(x)-v^\alpha(x')$ and of $\nabla v^\alpha$.

\begin{lemma}\label{main} There exists a constant
 $c(\ell,\hat c, \hat \eta)>0$  such that for all $x,\, x'\in H$
\begin{equation}\label{es1}
|v^\alpha(x)-v^\alpha(x')|\leq c (1+|x|^2+|x'|^2);
\end{equation}
and for all $x\in H$,
\begin{equation}\label{es2}
|\nabla v^\alpha(x)|\le c (1+|x|^2).
\end{equation}
\end{lemma}
We stress the fact that $c>0$ is independent of $\alpha$.

\proof
Set $$\tilde{\Upsilon}^\alpha(x)=\begin{cases} \dfrac{\psi(x,\nabla v^\alpha(x)G) -\psi(x,0)}{|\nabla v^\alpha(x)G|^2}
  \left(\nabla v^\alpha(x)G\right)^*, &  \hbox{ if $\; \nabla v^\alpha(x)G
 \neq 0$} \\ \\
  0,  \hbox{ if $\; \nabla v^\alpha(x)G
= 0.$}
 \end{cases}  $$
Then
$$\psi(X_t^x,Z_t^{x,\alpha})=\psi(X_t^x,0)+\tilde{\Upsilon}^\alpha(X_t^x)Z_t^{x,\alpha}.$$

From Proposition \ref{linear}, $\tilde{\Upsilon}^\alpha$ is the pointwise limit of a uniformly bounded sequence of Lipschitz  functions.

For all $T>0$, the couple of processes $(Y^{x,\alpha}, Z^{x,\alpha})$ is a solution to the following finite
horizon linear BSDE
\begin{equation}\label{BSDE_n_rew}
    \left\{
\begin{array}[c]{l}
-d Y^{x,\alpha}_t
={\psi}(X^x_t,0)dt+\tilde{\Upsilon}^{\alpha}(X^x_t)Z^{x,\alpha}_tdt
-\alpha Y^{x,\alpha}_tdt -Z^{x,\alpha}_t dW_t ,\text{ \ \ \
} t\in [0,T], \\
Y^{x,\alpha}_T  = v^{\alpha}(X^x_T).
\end{array}
\right.
\end{equation}
Since $\tilde{\Upsilon}^{\alpha}$ is bounded for all $T>0$, there exists
a unique probability $\hat{\mathbb{P}}^{x,\alpha,T}$ such that
$$\hat{W}^{x,\alpha}_t=\int_0^t \tilde{\Upsilon}^{\alpha}(X^x_s) ds +W_t $$ is
a $\hat{\mathbb{P}}^{x,\alpha,T}$-Wiener process for $t\in [0,T]$.
Consequently we have
$$ v^\alpha(x)= \hat{\mathbb{E}}^{x,\alpha,T} \left[ e^{-\alpha T}v^{\alpha}(X^x_T)
 +\int_0^T e^{-\alpha s} {\psi}(X^x_s,0)ds\right]
$$
where $\hat{\mathbb{E}}^{x,\alpha,T}$ denotes the expectation with
respect to $\hat{\mathbb{P}}^{x,\alpha,T}$.

Letting $T\rightarrow \infty$, as $|v^\alpha(x)|\le \frac{l}{\alpha}$, we get
$$ v^{\alpha}(x)=\lim_{T\rightarrow \infty} \hat{\mathbb{E}}^{x,\alpha,T} \left[
 \int_0^T e^{-\alpha s} {\psi}(X^x_s,0)ds\right].
$$
On the other hand, if we rewrite the forward equation (\ref{equazioneforwardmild}) with respect to
 $\hat{W}^{x,\alpha}$
 it turns out that $X^x$ verifies
 \begin{equation}
\left\{
\begin{array}[c]{l} d X^x_t  =A X^x_t
 d t+F(X^x_t) dt +G\tilde{\Upsilon}^{\alpha}(X^x_t)dt+G\hat{W}^{x,\alpha}_t , \\
\hat X^x_0  =x\in \, H.
\end{array}
\right.  \label{w-fsde}
\end{equation}
We denote by $\mathcal{P}^\alpha$ the associated Kolmogorov semigroup, i.e.,
$$\mathcal{P}_t^{\alpha}[\phi](x)=\hat{\mathbb{E}}^{x,\alpha,t}\phi(X_t^x).$$
Applying  Corollary \ref{cor-coupl} with $\Upsilon^\alpha=F+G\tilde{\Upsilon}^\alpha$ (which is also the pointwise limit of a sequence of Lipschitz functions), we obtain
$$|v^\alpha(x)-v^\alpha(x')|\leq \int_0^{\infty} e^{-\alpha t}
\left|\mathcal{P}^\alpha_t[\psi(\cdot,0)](x)-
 \mathcal{P}^\alpha_t[\psi(\cdot,0)](x')\right|dt\leq \frac{ \hat{c}l}{\hat{\eta}}
(1+|x^2|+|x'|^2) $$ where $\hat{c}$ and
$\hat{\eta}$ are independent of  $\alpha$. The proof of \eqref{es1} is now complete.

To prove \eqref{es2}, let us set
$$\bar{v}^\alpha(x)=v^\alpha(x)-v^\alpha(0).$$
Then, $\bar Y^{x,\alpha}_t = Y^{x,\alpha}_t -
Y^{0,\alpha}_0=\bar v ^{\alpha}(X^{x}_t)$ is the unique solution of
the finite horizon BSDE
$$
    \left\{
\begin{array}[c]{l}
-d \bar Y^{x,\alpha}_t  =\psi(X^x_t,Z^{x,\alpha}_t)dt -\alpha \bar{Y}_t^{x,\alpha}-\alpha v^\alpha(0)dt  -Z^{x,\alpha}_t dW_t
,\text{ \ \ \
} \\
Y^{x,\alpha}_1  = \bar v^\alpha(X^x_1).
\end{array}
\right.
$$
Note that in particular, in the above equation, $|\alpha v^\alpha(0)|\le l.$
By  Theorem 4.2 in \cite{FuTeBE}, $\bar{v}^\alpha$ is of class ${\cal G}^1$ and
there exists a constant $c(l,\hat{c},\hat{\eta})>0$ independent of $\alpha$ such that
$|\nabla v^\alpha(x)|\leq c(1+|x|^2)$, and the conclusion follows. \qed

\begin{remark}
As already mentioned in the introduction, the non degeneracy assumption on $G$ is essential in 
the proof of the gradient estimate on $v^\alpha$. More precisely, it is necessary to use the
Bismut-Elworthy formula from \cite{FuTeBE}.
\end{remark}

Now we are in position to state our main result in this section.

\begin{theorem} \label{main-EBSDE} Assume that the Hypotheses
\ref{ipotesiuno},
\ref{hypothesisroyer} and 3.2 hold. Moreover let $\bar \lambda$ be the
real number in (\ref{def-of-lambda}) below and define $\bar Y^x_t= \bar
v(X^x_t)$ (where $\overline{v}$ is a locally Lipschitz function with
$\overline{v}(0)=0$ defined in (\ref{def-of-v})). Then there
exists a process $\overline{Z}^{x}\in L_{\cal P, {\rm
loc}}^2(\Omega;L^2(0,\infty;\Xi^*))$
  such that $\mathbb{P}$-a.s.  the EBSDE
 (\ref{EBSDE}) is satisfied by
 $(\bar Y^x,\bar Z^x, \bar \lambda)$ for all $0\leq t\leq T$.

Moreover   $\overline{v}$ is of class ${\cal G}^1$, $|\nabla \overline{v}(x)|\le c(1+|x|^2)$, and
$\overline{Z}^{x}_t=\nabla\overline{v} (X^x_t)G$.
\end{theorem}

\proof  Let us set again $\bar{v}^\alpha(x)=v^\alpha(x)-v^\alpha(0).$  By  Lemma \ref{lemmaroyer} and relation (\ref{es1})
      we can construct, by a diagonal procedure, a sequence $\alpha_n\searrow 0$ such that for all $x$ in a
countable dense subset $D\subset H$
     \begin{equation}\label{def-of-lambda}
    {\overline{v}}^{\alpha_n}(x)\rightarrow \overline{v}(x),\qquad
\alpha_n v^{\alpha_n}(0)\rightarrow \overline{\lambda},
    \end{equation}
for a suitable  function $  \overline{v}: D \rightarrow
\mathbb{R}$ and for a suitable real number $\overline{\lambda}$.

 Moreover, by Lemma \ref{main}, $ | \overline{v}^{\alpha}(x)-
\overline{v}^{\alpha}(x')|\leq    c(1+|x|^2+|x'|^2)|x-x'|$   for all
$x,x'\in H$ and all $\alpha>0$. So $\overline{v}$ can be extended
to a locally Lipschitz function defined on the whole $H$ with  $| \overline{v}(x)-
\overline{v}(x')|\leq    c(1+|x|^2+|x'|^2)|x-x'|$ and
\begin{equation}\label{def-of-v} {\overline{v}}^{\alpha_n}(x)\rightarrow
\overline{v}(x),\qquad  x\in H.\end{equation}

Clearly we have,
$\mathbb{P}$-a.s.,
\begin{equation}\label{equation-proof-main-1}
 \overline{Y}^{x,\alpha}_t=\overline{Y}^{x,\alpha}_T +\int_t^T(\psi(X^{x}_\sigma,
Z^{x,\alpha}_\sigma)-\alpha \overline{Y}^{x,\alpha}_\sigma-\alpha
{v}^{\alpha}(0))d\sigma -\int_t^T Z^{x,\alpha}_\sigma dW_\sigma,
\; 0\le t\le T <\infty.
\end{equation}
Since $|\bar v^{\alpha}(x)|\le c(1+|x|^2)$, inequality
(\ref{stimadeimomentidix})   ensures that
$\mathbb{E}\sup_{t\in[0,T]}\left[\sup_{\alpha>0}
|\overline{Y}^{x,\alpha}_t|^2\right]< +\infty$  for any $T>0$.
Thus, if we define $\overline{Y}^x=\overline{v}(X^x)$, then by
dominated convergence theorem
$$\mathbb{E} \int_0^T |\overline{Y}^{x,\alpha_n}_t -\overline{Y}^{x}_t|^2 dt
 \rightarrow 0\quad \hbox{and}\quad
\mathbb{E} |\overline{Y}^{x,\alpha_n}_T-\overline{Y}^{x}_T|^2
\rightarrow 0
$$
as $n\rightarrow \infty$ (where $\alpha_n \searrow 0$ is a
sequence for which (\ref{def-of-lambda}) and (\ref{def-of-v})
hold).

We claim now that there exists $\overline{Z}^{x}\in L_{\cal P,
{\rm loc}}^2(\Omega;L^2(0,\infty;\Xi^*))$ such that
 $$\mathbb{E} \int_0^T |{Z}^{x,\alpha_n}_t -\overline{Z}^{x}_t|_{\Xi^*}^2 dt
  \rightarrow 0.$$
Let  $\tilde{Y}={\bar Y}^{x,\alpha_n}-{\bar Y}^{x,\alpha_m}$,
$\tilde{Z}={Z}^{x,\alpha_n}-{Z}^{x,\alpha_m}$. Applying It\^o's
rule to $\tilde{Y}^2$ we get, by standard computations, that
$$\tilde{Y}^2_0+\mathbb{E}\int_0^T |\tilde{Z}_t|_{\Xi^*}^2 dt
=\mathbb{E}{\tilde Y}^2_T + 2\mathbb{E}\int_0^T \tilde \psi_t
\tilde Y_t dt -2 \mathbb{E}\int_0^T \left[\alpha_n
{Y}^{x,\alpha_n}_t - \alpha_m {Y}^{x,\alpha_m}_t\right] \tilde
Y_t\,dt,
$$
where  $\tilde
\psi_t=\psi(X^x_t,Z^{x,\alpha_n}_t)-\psi(X^x_t,Z^{x,\alpha_m}_t)$.
We notice that $|\tilde\psi_t| \leq  l|\tilde Z _t|$ and
$\alpha_n |{Y}^{x,\alpha_n}_t|$ $\leq l$. Thus
$$
\mathbb{E}\int_0^T |\tilde{Z}_t|_{\Xi^*}^2 dt \leq c\left[
\mathbb{E} (\tilde Y_T)^2 +\mathbb{E}\int_0^T (\tilde{Y}_t)^2
dt +\mathbb{E}\int_0^T |\tilde{Y}_t| dt \right].$$ It follows
that the sequence $\{{Z}^{x,\alpha_n}\}$ is Cauchy in
$L^2(\Omega;L^2(0,T;\Xi^*))$ for all $T>0$ and our claim is
proved.

Now we can  pass to the limit as $n\rightarrow \infty$ in equation
(\ref{equation-proof-main-1}) to obtain
\begin{equation}\label{equation-proof-main-2}
 \overline{Y}^{x}_t=\overline{Y}^{x}_T +\int_t^T(\psi(X^{x}_\sigma,
\overline{Z}^{x}_\sigma)-\overline{\lambda })d\sigma-\int_t^T
\overline{Z}^{x}_\sigma dW_\sigma, \quad 0\le t\le T <\infty.
\end{equation}
We notice  that the above equation also ensures continuity of the
trajectories of $\overline{Y}$.

Finally, the couple of processes $(\overline{Y}^x,\overline{Z}^x)$  is the unique solution of
the finite horizon BSDE
$$
    \left\{
\begin{array}[c]{l}
-d \bar Y^{x}_t  =(\psi(X^x_t,\overline{Z}^x_t)-\overline{\lambda})dt  - \bar Z^{x}_t dW_t
,\text{ \ \ \
} \\
\bar Y^{x}_1  = \bar v(X^x_1).
\end{array}
\right.
$$

Once again, by  Theorem 4.2 in \cite{FuTeBE}, we conclude the proof. \qed

\begin{remark}\begin{em} The solution we
have constructed above has the following ``quadratic growth''
property with respect to $X$: there exists $c>0$ such that,
$\mathbb{P}$-a.s.,
\begin{equation}\label{growt-of-Y}
|\overline{Y}^x_t|\leq c (1+ |X^x_t|^2),  \hbox{ for all $t\geq 0$}.
\end{equation}
\end{em}
\end{remark}
If we require similar conditions then we immediately obtain
uniqueness of $\lambda$.
\begin{theorem}\label{th-uniq-lambda} Assume that the Hypotheses \ref{ipotesiuno}, \ref{hypothesisroyer} and 3.2
hold true.  Moreover suppose that, for
some $x\in H$, the triple $(Y',Z',\lambda')$ verifies
$\mathbb{P}$-a.s. equation
 (\ref{EBSDE}) for all $0\leq t\leq T$,
where
 $Y'$ is a progressively measurable continuous process, $Z'$ is a process
 in $L_{\cal P, {\rm loc}}^2(\Omega;L^2(0,\infty;\Xi^*))$ and
 $\lambda'\in \mathbb{R}$.
 Finally assume that there exists $c_x>0$ (that may depend
 on $x$) such that for some $p\ge 1$,
$\mathbb{P}$-a.s.
$$
 |Y'_t|\leq c_x (1+|X^x_t|^p) , \hbox{ for all $t\geq 0$}.
$$ Then $\lambda'=\bar \lambda$.
\end{theorem}

The proof of the above theorem is similar to that of Theorem 4.6 in \cite{FHT-Erg}, so we omit it here.

$ $

The solution obtained in Theorem \ref{main-EBSDE} has moreover the property that
processes $Y^x$ and $Z^x$ are deterministic functions of $X^x$. We refer to such solutions as to
``Markovian" solution of the EBSDEs.

We prove that the Markovian solution is unique.

\begin{theorem}\label{markov} Let $(v,\zeta)$, $(\tilde v,\tilde \zeta)$ two
couples of functions with $v, \tilde v :H \rightarrow \mathbb{R}$,
continuous, with $ |v(x)|\leq c (1+|x|^2)$, $ |\tilde v(x)|\leq c
(1+|x|^2)$, $v(0)=\tilde v(0)=0$ and $\zeta,\,\tilde \zeta$ continuous
from $H$ to $\Xi^*$ endowed with the weak$^*$ topology verifying $
|\zeta(x)|\leq c (1+|x|^2)$, $ |\tilde \zeta(x)|\leq c (1+|x|^2)$.

Assume that for some constants $\lambda$, $\tilde \lambda$ and all
$x\in H$,
$(v(X^x_t),\zeta(X^x_t), \lambda)$, $
(\tilde v(X^x_t),\tilde \zeta(X^x_t), \tilde \lambda) $ verify the EBSDE (\ref{EBSDE}),  then
$\lambda= \tilde \lambda,\; v= \tilde v,\; \zeta=\tilde\zeta.$
\end{theorem}
\proof The equality $\lambda= \tilde \lambda$ comes from Theorem \ref{th-uniq-lambda}.

Then let $\bar {Y}^x_t=v(X^x_t)-\tilde v(X^x_t)$, $\bar
{Z}^x_t=\zeta(X^x_t)-\tilde \zeta(X^x_t)$ and $\tilde\Upsilon$ be
defined in Proposition \ref{linear}. We have
$$-d\bar {Y}^x_t=\tilde \Upsilon(X^x_t)\bar {Z}^x_t dt - \bar {Z}^x_t dW_t=
-\bar {Z}^x_t dW'_t$$ where $W'_t=-\int_0^t \Upsilon(X^x_s) ds+W_t$ is a
Wiener process in [0,T] under the probability $\bar {\mathbb{P}}^{x,T}$.

Moreover, under $\bar {\mathbb{P}}^{x,T}$, $X^x$ satisfies
equation (\ref{sde-coup}), in $[0,T]$, with, as before $\Upsilon=G
\widetilde{\Upsilon}+F$. Thus, from Proposition \ref{esisteunicox}, it holds that  for all $p\geq 1$, and all $x\in H$
$$\bar{\mathbb{E}}^{x,T}|X^x_t|^p\leq c(1+|x|^p), \forall 0\leq
t\leq T,$$ where $c>0$ depends on $p,\gamma,M$ and $l|G|+\sup_x |F(x)|$, and is independent of $T$.
 Thus the growth conditions on $\zeta$ and $\tilde \zeta$
implies that, for all $T>0$, $\bar{\mathbb{E}}^{x,T}\int_0^T |\bar
{Z}^x_t|^2 dt <\infty$.

Let $\tau =\inf\{t : | X^x_t|< \epsilon\}$ then for all $T>0$
$$ \bar {Y}^x_0= \bar{\mathbb{E}}^{x,T} \bar {Y}^x_{T\wedge \tau}.$$
For any $\delta>0$, there exists $\epsilon>0$ such that $|v(x)-\tilde v(x)|\leq \delta$ if
$|x|\leq \epsilon$. Then for a constant $c>0$,
\begin{eqnarray*}
|\bar Y_0^x|=|\bar{\mathbb{E}}^{x,T} \bar {Y}^x_{T\wedge \tau}|&\leq&\bar{\mathbb{E}}^{x,T} |\bar {Y}^x_{\tau}|1_{\{\tau<T\}}+\bar{\mathbb{E}}^{x,T} |\bar {Y}^x_{T}|1_{\{\tau\ge T\}}\\
&\le& \delta +
\left(\bar{\mathbb{P}}^{x,T}\{ \tau\geq T\}\right)^{1/2}
\left(\bar{\mathbb{E}}^{x,T}\{  |\bar {Y}^x_{T}|^2\}\right)^{1/2}\\
&\le&
 \delta +
\left(\bar{\mathbb{P}}^{x,T}\{ \tau\geq T\}\right)^{1/2}
\left(\bar{\mathbb{E}}^{x,T}\{ 1+|X^x_T|^4\}\right)^{1/2}.
\end{eqnarray*}
Noting that, by Theorem \ref{recurrence}, $\lim_{T\rightarrow \infty}\bar{\mathbb{P}}^{x,T}\{ \tau\geq T\}=0$ and sending $T$ to $\infty$ in the last inequality, we obtain that $|\bar {Y}^x_0|\leq \delta$ and the claim follows from the arbitrarity of $\delta$. \qed

\section{Ergodic HJB equations}\label{sec-HJB}
We briefly show here that as $\bar v(x)=\bar Y_0^x$ in Theorem \ref{main-EBSDE}  is of class
${\cal G}^1$, the couple $(\bar v,\bar \lambda)$  is a mild solution of
the following ``ergodic''  Hamilton-Jacobi-Bellman equation:
\begin{equation}
\mathcal{L}v(x)
+\psi\left( x,\nabla v(x)  G\right)  = \lambda, \quad x\in H,  \label{hjb}%
\end{equation}
where  the linear operator $\mathcal{L}$ is formally defined by
\[
\mathcal{L}f\left(  x\right)  =\frac{1}{2}Trace\left(
GG^{\ast}\nabla ^{2}f\left(  x\right)  \right)  +\langle Ax,\nabla
f\left(  x\right) \rangle+\langle F\left(  x\right)
,\nabla f\left(  x\right) \rangle.
\]
We notice that we can define the transition semigroup
 $(P_t)_{t\geq 0}$ corresponding to $X$ by the formula $P_t[\phi](x)=E\phi(X_t^x)$
for all measurable functions $\phi:E\to\mathbb{ R}$ having
polynomial growth, and we notice that $\mathcal{L}$ is the formal
generator of $(P_t)_{t\geq 0}$.

 Since we are dealing with an elliptic equation it is natural to consider
$(v,\lambda)$ as a mild solution of equation (\ref{hjb}) if and
only if, for arbitrary $T>0$, $v(x)$ coincides with the  mild
 solution $u(t,x)$ of the corresponding parabolic equation
 having $v$ as a terminal condition:
\begin{equation}\left\{
\begin{array}{l}
 \frac{\partial u(t,x)}{\partial t}+\mathcal{L}u\left(  t,x\right)
+\psi\left(  x,\nabla u\left(  t,x\right)  G\right)
 -\lambda=0, \quad t\in [0,T],\; x\in H,  \\ \\
u(T,x)=v(x), \quad  x\in H.
 \end{array}\right. \label{hjb-parab}
\end{equation}
Thus we are led to the following definition:
\begin{definition}
\label{defsolmildkolmo} A pair $(v,\lambda)$ ($v: H\rightarrow
\mathbb{R}$ and $\lambda\in \mathbb{R}$) is a mild solution of the
Hamilton-Jacobi-Bellman equation (\ref{hjb}) if the following are
satisfied:

\begin{enumerate}
\item $v\in\mathcal{G}^{1}\left(  H,\R \right)  $;

\item  there exists $C>0$ such that $\left|  \nabla v\left(
x\right)
\right|  \leq C\left(  1+\left|  x\right|^{p}\right)  $ for every  $x\in H$ and some
$p\ge 1$;

\item for $0\le t\le T$  and $x\in H$,
\begin{equation}
v(x)=P_{T-t}\left[  v\right]  \left(  x\right)
+\int_{t}^{T}\left(P_{s -t }\left[  \psi(\cdot,\nabla v\left(
\cdot\right)  G)\right]  \left( x\right) -\lambda \right)   \,ds.
\label{mild sol hjb}
\end{equation}

\end{enumerate}
\end{definition}

Theorems \ref{main-EBSDE} and \ref{markov} immediately yield
existence and uniqueness of the mild solution of equation (\ref{hjb}).

\begin{theorem}\label{th-EHJB}
Assume that Hypotheses \ref{ipotesiuno} and
\ref{hypothesisroyer} hold.

Then $(\bar v,  \bar\lambda)$ is a   mild solution of the
Hamilton-Jacobi-Bellman equation (\ref{hjb}).

Conversely, if $(v,\lambda)$ is a   mild solution of
 (\ref{hjb}) then, setting $ Y^x_t=
v(X^x_t)$ and ${Z}^{x}_t=
\nabla v( X^x_t)  G$,
the triple
 $( Y^x, Z^x,  \lambda)$ is a solution of
  the EBSDE
 (\ref{EBSDE}), which implies the uniqueness of mild solution in the sense that if
$(\bar v^{i},  \bar\lambda)$, $i=1,2$ are   mild solutions of the
Hamilton-Jacobi-Bellman equation (\ref{hjb}) then
$v^1(X^x_t)=v^2(X^x_t)$  and $\nabla v^1( X^x_t)  G=\nabla v^2( X^x_t)  G$  $\mathbb{P}$- a.s. for a.e. $t\geq 0$.
\end{theorem}
\proof: The proof is identical to the one of Theorem 6.2 in \cite{FHT-Erg}.

\section{Ergodic control}\label{sec-contr}
We fix  a bounded function $F:H\rightarrow H$ and denote 
 by ${X}^{x}$  the solution of equation (\ref{equazioneforwardmild})
 with $\Upsilon=F$.

Assume that the Hypotheses \ref{ipotesiuno} and 3.2 hold.
 Let $U$ be a separable
 metric space. We define a control $u$ as an
$(\calf_t)$-progressively measurable $U$-valued process.  The cost
 corresponding to a given control
is defined in the following way. We assume that the functions
$R:U\rightarrow \Xi^*$ and $L:H\times U \rightarrow \R$ are
measurable and satisfy, for some constant $c>0$,
\begin{equation}\label{condcosto}
|R(u)|\leq c,\quad |L(x,u)|\leq c, \quad |L(x,u)-L(x',u)|\leq
c\,|x-x'|,\qquad u\in U,\,x,x'\in H.
\end{equation}
Given an arbitrary control $u$ and $T>0$, we introduce the
Girsanov density
$$ \rho_T^u=\exp\left(\int_0^T R(u_s)dW_s
-\frac{1}{2}\int_0^T |R(u_s)|_{\Xi^*}^2 ds\right)$$ and the
probability $\mathbb P_T^u=\rho_T^u\mathbb P$ on $\calf_T$. The
ergodic cost  corresponding to $u$ and the starting point $x\in H$
is
\begin{equation}\label{def-ergodic-cost}
  J(x,u)=\limsup_{T\rightarrow\infty}\frac{1}{T} \mathbb
E^{u,T}\int_0^T L(X_s^x,u_s)ds,
\end{equation}
where $\mathbb E^{u,T}$  denotes expectation with respect to
$\mathbb P_T^u$. We notice that $W_t^u=W_t-\int_0^t R(u_s)ds$ is a
Wiener process on $[0,T]$ under $\mathbb P^u_T$ and that
$$dX_t^x=(AX_t^x+F(X_t^x))dt+G(dW_t^u+R(u_t)dt),
\quad t\in [0,T]$$ and this justifies our formulation of the
control problem. Our purpose is to minimize the cost over all
controls.

 To this purpose we first define the Hamiltonian in the
usual way
\begin{equation}\label{defhamiton}
\psi(x,z)=\inf_{u\in U}\{L(x,u)+z R(u)\},\qquad x\in H,\,z\in
\Xi^*,
\end{equation}
and we remark that  if, for all $ x,z$, the  infimum is attained
in (\ref{defhamiton}) then  by the Filippov Theorem, see \cite{McS-War}, there exists a measurable function
$\gamma:H\times \Xi^*\rightarrow U$ such that
$$\psi(x,z)=l(x,\gamma(x,z))+z R(\gamma(x,z)).$$

We notice that under  the present assumptions $\psi$ is a
Lipschitz function and $\psi(\cdot,0)$ is bounded (here the fact
that $R$ depends only on $u$ is used). So if we assume the Hypotheses
\ref{ipotesiuno} and 3.2, then in Theorem
\ref{main-EBSDE} we have constructed, for every $x\in H$, a triple
\begin{equation}\label{richiamoebsde}
(\bar Y^x,\bar Z^x, \bar \lambda)= (\bar v (X^x),\bar \zeta(X^x),
\bar \lambda)
\end{equation} solution to
 the EBSDE
 (\ref{EBSDE}).

\begin{theorem}\label{Th-main-control}
Assume that the Hypotheses \ref{ipotesiuno} and 3.2 hold, and that
(\ref{condcosto}) holds as well.

Moreover suppose that, for some $x\in H$, a triple $(Y,Z,\lambda)$
verifies $\mathbb{P}$-a.s. equation
 (\ref{EBSDE}) for all $0\leq t\leq T$,
where
 $Y$ is a progressively measurable continuous process, $Z$ is a process
 in $L_{\cal P, {\rm loc}}^2(\Omega;L^2(0,\infty;\Xi^*))$ and
 $\lambda\in \mathbb{R}$.
 Finally assume that there exists $c_x>0$ (that may depend
 on $x$) such that
$\mathbb{P}$-a.s.
$$
 |Y_t|\leq c_x (1+|X^x_t|^2) , \hbox{ for all $t\geq 0$}.
$$

Then the following holds:
\begin{enumerate}
 \item[(i)] For arbitrary control
 $u$ we have $J(x,u)\ge \lambda=\bar\lambda,$
and the equality holds if  $L(X_t^x,u_t)+Z_t
R(u_t)=\psi(X_t^x,Z_t)$, $\P$-a.s. for almost every $t$.

\item[(ii)] If the  infimum is attained in (\ref{defhamiton}) then
the control $\bar u_t=\gamma(X_t^x,Z_t)$ verifies $J(x,\bar u)=
\bar\lambda.$
\end{enumerate}

In particular, for the solution (\ref{richiamoebsde}) mentioned
above, we have:
\begin{enumerate}
 \item[(iii)] For arbitrary control
 $u$ we have $J(x,u)=\bar\lambda$ if
$L(X_t^x,u_t)+\bar\zeta (X_t^x) R(u_t)=\psi(X_t^x,\bar \zeta
(X_t^x))$, $\P$-a.s. for almost every $t$. \item[(iv)] If the
infimum is attained in (\ref{defhamiton}) then the control $\bar
u_t=\gamma(X_t^x,\bar\zeta (X_t^x))$ verifies $J(x,\bar u)=
\bar\lambda.$
\end{enumerate}
\end{theorem}
\proof: The proof is identical to the one of Theorem 7.1 in \cite{FHT-Erg}.

\begin{example}\label{ex-heat-eq} We consider here an ergodic optimal control when the state
equation is a stochastic heat equation. The difference with respect to the same
 example in \cite{FHT-Erg} is that, if we have non-degenerate noise,
  we do not need to assume that the non-linearity $f$ is decreasing.
  Namely we consider the following state equation
\begin{equation}
\left\{
\begin{array}
[c]{l} d_{t }X^{u}\left(  t ,\xi\right)  =\left[
\frac{\partial^{2}}{\partial \xi^{2}}X^{u}\left(  t ,\xi\right)
+f\left(  \xi,X^{u}\left(  t  ,\xi\right)  \right)
+r(\xi,u\left( t ,\xi\right) ) \right] dt+
\sigma(\xi) \dot{W}\left(
t ,\xi\right)  dt ,\\
X^{u}\left(  t ,0\right)  =X^{u}\left(  t ,1\right)  =0,\\
X^{u}\left(  t,\xi\right)  =x_{0}\left(  \xi\right)  ,
\end{array}
\right.  \label{heat equation}
\end{equation}
where $\dot{W}\left( t ,\xi\right)  $ is a
space-time white noise on $[0,\infty[ \times\left[
0,1\right] $.

\noindent We also introduce the cost functional
\begin{equation}
J\left(  x,u\right)  = \limsup_{T\rightarrow\infty}\dfrac{1}{T}
\mathbb{E}\int_{0}^{T}\int_{0}^{1}l\left(  \xi ,X^{u}_s\left(
\xi\right)  ,u_s(\xi)\right)  d\xi  \, ds.
 \label{heat costo diri}
\end{equation}
An admissible control $u\left(  t  ,\xi\right)  $ is a predictable
process $u : \Omega\times [0,\infty[\times[0,1]$.

$ $

Then (the reduction to the abstract infinite dimensional framework is as in \cite{FuTeBE} Section 5.1)
 the results of Theorem \ref{Th-main-control}
 can be applied under the following assumptions:

\begin{enumerate}
    \item  $f:[0,1]\times \R \rightarrow \R$ is measurable and bounded and
    $$| f(\xi,\eta_{1})- f(\xi,\eta_{2})|\leq
    c_f|\eta_{2}-\eta_{1}|,
$$
    for a  suitable constant $c_f$, almost all
$   \xi\in [0,1]$,  and all
    $\eta_{1}, \eta_{2}\in \mathbb{R}$.
    Moreover we assume that $ f(\xi,\cdot)\in C^{1}(\R)$ for a.a.  $\xi\in [0,1]$.
    \item $\sigma: [0,1]\rightarrow \R$ is measurable and bounded. Moreover
    $
c_{\sigma}\leq  | \sigma(\xi)|$,
 for a.a.  $\xi\in [0,1]$ and a suitable constant $c_{\sigma}>0$.
    \item  $r: [0,1]\times \R$ is measurable and bounded
 and,  for a.a.  $\xi\in [0,1]$, the map $r(\xi,\cdot): \R\to\R$ is continuous.
    \item  $l:[0,1]\times \R^2 $ is measurable and bounded and, for a.a.  $\xi\in [0,1]$, the
    map  $l(\xi,\cdot,\cdot):\R^2 \to \R$ is
continuous.
\item  $x_0\in L^{2}([0,1])$.
\end{enumerate}

\end{example}

\section{Appendix}\label{appendix}

\subsection{Proof of the coupling estimates}\label{sec-app-coupl}
\proof
We use a coupling argument. The precise result needed here cannot be found in the literature and
we give the full proof for completeness. In particular, we show that the constant $\hat c$ and
$\hat \eta$ in \eqref{coupling-estimate} depend only on the supremum of $\Upsilon$. Our proof
is a mix of arguments found in \cite{EMS}, \cite{kuksin-shirikyan}, \cite{mattingly}, 
\cite{shirikyan} and \cite{odasso}. Note that the more analytical method from \cite{DPDT05} could also be used.

To prove that the laws of the solutions starting from different initial data get closer in variation, we first
wait that these solutions enter a fixed ball. Then, we construct a coupling of solutions starting from
initial data in this ball. Iterating this argument we obtain the result.

In this section, $\kappa_i,\; i=0,1\dots$ denotes a constant which depends only on
$\gamma,M,G$ and $L_0=\sup_{x\in H}|F(x)|$.

{\bf Step 1:}
Let $Z_t=\int_0^t e^{A(t-s)}GdW(s)$ and $\rho=X^x-Z$. Thanks to Hypothesis \ref{ipotesiuno}, we obtain by taking the scalar product of the equation verified by $\rho$ with $\rho$:
$$
\frac12\frac{d}{dt} |\rho|^2 +k |\rho|^2\le L|\rho| \le \frac{L_0^2}{2k} + \frac{k}2 |\rho|^2
$$
with $L_0=\sup_{x\in H}|\Upsilon(x)|$. By Gronwall's lemma
$$
|\rho_t|^2\le e^{-kt}(|x|^2+\frac{L_0^2}{k^2}).
$$
Moreover
$$
\E(|Z_t|^2)=\int_0^t |e^{A(t-s)}G|_{L_2(H;H)}^2ds \le |G|_{{\mathcal L}(H)}^2 \int_0^t |e^{As}|_{L_2(H;H)}^2ds.
$$
It follows
\begin{equation}\label{Ap1}
\mathbb{E}|X^x_t|^2\le 2(|x|^2e^{-kt}+\kappa_1),
\end{equation}
with a constant $\kappa_1$ depending only on $\gamma,M,G$ and $L_0=\sup_{x\in H}|F(x)|$ but independent of $t>0$.
By the Markov property:
\begin{equation}\label{e7.2}
\E(|X^x_{(k+1)T}|^2|{\mathcal F}_{kT})\le  2(|X^x_{kT}|^2e^{-kT})+\kappa_1,\; k\ge 0.
\end{equation}
Let us define for $R\ge 0$
$$
C_k=\{ |X^x_{kT}|^2\ge R\},\quad B_k=\cap_{j=0}^k C_j.
$$
By Chebychev's inequality
\begin{equation}\label{e7.3}
\P(C_{k+1}|{\mathcal F}_{kT})\le \frac{2e^{-kT}}{R}|X^x_{kT}|^2+\frac{\kappa_1}{R}.
\end{equation}
Multiply \eqref{e7.2} and \eqref{e7.3} by ${\bf 1}_{B_k}$ and take the expectation to obtain,
since ${\bf 1}_{B_{k+1}}\le {\bf 1}_{B_k}$,
$$
\left(
\begin{array}{c}
\E(|X^x_{(k+1)T}|^2{\bf 1}_{B_{k+1}})\\
\\
\P(B_{k+1})
\end{array}
\right)
\le A
\left(
\begin{array}{c}
\E(|X^x_{kT}|^2{\bf 1}_{B_{k}})\\
\\
\P(B_{k})
\end{array}
\right)
$$
where
$$
A=\left(
\begin{array}{cc}
2e^{-kT} & \quad \kappa_1\\
\\
\displaystyle \frac{2e^{-kT}}{R} & \quad \displaystyle \frac{\kappa_1}{R}
\end{array}
\right).
$$
Let $R=4\kappa_1$ and choose $T$ such that $e^{-kT}= \frac{1}8$. The eigenvalues of $A$ are $0$ and
$2 e^{-kT}+\frac14\le \frac18$. We deduce that
$$
\P(B_k)\le \kappa_2 \left(\frac18\right)^k (1+ |x|^2),
$$
with $\kappa_2$ depending only on $\kappa_1$.
Defining
$$\tau=\inf\{kT; |X^x_{kT}|^2\le R\},$$
it follows
$$\mathbb{P}(\tau\ge kT)\le \P(B_k)\le \kappa_2 \left(\frac18\right)^k (1+ |x|^2).$$
Thus, for $\eta T<2\ln 2$,
\begin{equation}\label{Ap2}
\mathbb E(e^{\eta\tau})\le \kappa_3(1+|x|^2).
\end{equation}

{\bf Step 2:}
We construct a coupling for $x,y\in B_R$, the ball of center $0$ and radius $R$, $x\not=y$. Let $T\ge 0$ to be chosen below, we denote by $\mu_1$ the law of $X^x$ and $\mu_2$
the law of $X^y$ on $[0,T]$. Set
$\tilde{X}_t=X_t^y+\frac{T-t}{T}e^{At}(x-y)$ and denote by $\tilde{\mu_2}$ the law of $\tilde{X}$
on $[0,T]$.
Then
$$d\tilde{X}=(A\tilde{X}+F(\tilde{X}))dt+d\tilde{W},$$
where $\tilde{W}_t=W_t-\int_0^t h(s)ds$, with $h(s)=F(\tilde{X}_s)-F(X_s^y)-\frac{1}{T}e^{At}(x-y)$.
By Girsanov's formula, $\tilde{W}$ is a Wiener process under a new probability measure $\tilde{\mathbb P}$.
Therefore, under $\tilde{\mathbb P}$, $\tilde{X}$ has the law $\mu_1$ while under $\mathbb P$ it has the law $\tilde{\mu}_2$.
Of course $\mu_1$ and $\tilde{\mu_2}$ are equivalent. Since $|h(t)|\le 2L_0+2\frac{R}{T}$, we deduce that
$$\int_H \left(\frac{d\tilde{\mu}_2}{d\mu_1}\right)^{3}d\mu_1\le \kappa_4.$$

We need the following result (see for instance \cite{mattingly}).
\begin{proposition}\label{p5.1}
Let $(\mu_1,\mu_2)$ be two probability measures on a Banach space $E$ then
$$
\|\mu_1-\mu_2\|_{TV}=\min \P(Z_1\ne Z_2)
$$
where the minimum is taken on all coupling $(Z_1,Z_2)$ of $(\mu_1,\mu_2)$. Moreover, there exists a
coupling which realizes the infimum. We say that it is a maximal coupling. It satisfies\footnote{Recall that if $\mu_1,\mu_2$ are absolutely continuous with respect to a measure
$\mu$ (for instance $\mu=\mu_1+\mu_2$), we have
$$
\begin{array}{l}
d(\mu_1\wedge \mu_2)=(\frac{d\mu_1}{d\mu}\wedge \frac{d\mu_2}{d\mu})d\mu.
\end{array}
$$}
$$
\P(Z_1=Z_2, \; Z_1\in \Gamma)=\mu_1\wedge \mu_2 (\Gamma),\; \Gamma\in \mathcal B(E).
$$
Moreover, if $\mu_1$ and $\mu_2$ are equivalent and   if
$$
\int_E \left(\frac{d\mu_2}{d\mu_1}\right)^{p+1}d\mu_1 \le C
$$
for some $p>1$
and $C>1$
then
$$
\P(Z_1=Z_2)=\mu_1\wedge \mu_2(E)\ge (1-\frac1p)\left(\frac{1}{pC}\right)^{1/(p-1)}.
$$
\end{proposition}
We deduce the existence of a coupling $(V^{1,x,y},\tilde V^{2,x,y})$ of $(\mu_1,\tilde \mu_2)$ such that
$$
\P(V^{1,x,y}=\tilde V^{2,x,y})\ge \frac{1}{4\kappa_4}.
$$
Clearly, $(V^{1,x,y}_t,V^{2,x,y}_t=\tilde V^{2,x,y}_t-\frac{T-t}Te^{At}(x-y))_{t\in [0,T]}$ is a coupling of
$(\mu_1,\mu_2)$ on $[0,T]$ and
\begin{equation}\label{Ap3}
\P(V^{1,x,y}_T= V^{2,x,y}_T)\ge \P(V^{1,x,y}= \tilde V^{2,x,y})\ge\frac{1}{4\kappa_4}.
\end{equation}
{\bf Step 3:} We now construct a coupling for any initial data.
For $x=y$, we set
$$
(V^{1,x,x}_t,V^{2,x,x}_t) = (X^x_t,X^x_t),\; t\in [0,T].
$$
If $x$ or $y$ is not in $B_R$, we set
$$
(V^{1,x,y}_t,V^{2,x,y}_t)=(X^x_t,\bar X^y_t),\; t\in [0,T],
$$
where
$\bar X^y$ is the solution of equation \eqref{sde-coup} driven by a Wiener process
$\bar W$ independent of $W$.
The coupling of the laws of $X^x$, $X^y$ for $t\ge 0$ is defined recursively
by the formula
$$V^{i,x,y}_{nT+t}=V_t^{i,V^{1,x,y}_{nT},V^{1,x,y}_{nT}}, \; t\in [0,T],\quad i=1,2.$$
We then define the following sequence of stopping times:
$$L_m=\inf\{l>L_{m-1}, V^{1,x,y}_{lT}, V^{2,x,y}_{lT}\in B_R\}$$
with $L_0=0$. Evidently, (\ref{Ap2}) can be generalized to two solutions and we have:
$$\mathbb E(e^{\eta L_1T})\le \kappa_3(1+|x|^2+|y|^2)$$
and
$$
\mathbb E(e^{\eta(L_{m+1}-L_m)T}|{\cal F}_{L_mT})\le \kappa_3(1+|V^{1,x,y}_{L_mT}|^2+|V^{2,x,y}_{L_mT}|^2).
$$
It follows
$$
\begin{array}{ll}
\E(e^{\eta L_{m+1}T} )&\displaystyle\le \kappa_3\E \left(e^{\eta L_m T} (1+|V^{1,x,y}_{L_mT}|^2+|V^{2,x,y}_{L_mT}|^2)\right)\\
&\le  \displaystyle \kappa_3(1+2R^2)\E \left(e^{\eta L_m T} \right)
\end{array}
$$
and
$$
\E \left(e^{\eta L_m T} \right) \le \kappa_3^l(1+2R^2)^{l-1}(1+|x|^2+|y|^2).
$$
Set now
$$
\ell_0=\inf\{ l,\; V^{1,x,y}_{(L_l+1)T}=V^{2,x,y}_{(L_l+1)T}  \}.
$$
Since $ V^{1,x,y}_{L_lT}=V^{2,x,y}_{L_lT}\in B_R$, we have by \eqref{Ap3}
$$
\P(\ell_0>l+1\big| \ell_0>l)\le(1-\frac1{4\kappa_4}).
$$
Since $\P(\ell_0>l+1)=\P(\ell_0>l+1\big| \ell_0>l)\P(\ell_0>l)$,
we obtain
$$
\P(\ell_0>l )\le(1-\frac1{4\kappa_4})^l.
$$
Then for $\gamma\ge 0$
$$
\begin{array}{ll}
\E(e^{\gamma L_{\ell_0}T})&\le \displaystyle\sum_{l\ge 0}\E(e^{\gamma L_lT} {\bf 1}_{l=\ell_0})\\
&\displaystyle \le \sum_{l\ge 0} \P(l=\ell_0)^{1-\gamma/\eta}(\E(e^{\eta L_lT}))^{\gamma/\eta}\\
&\le\displaystyle\sum_{l\ge 0} \left(1-\frac1{4\kappa_4}\right)^{(l-1)(1-\gamma/\eta)}
\left[   \kappa_3^l(1+2R^2)^{l-1}(1+|x|^2+|y|^2) \right]^{\gamma/\eta }.
\end{array}
$$
We choose $\gamma\le \eta $ such that
$$
\left(1-\frac1{4\kappa_4}\right)^{1-\gamma/\eta}
\left[   \kappa-3(1+2R^2) \right]^{\gamma/\eta }<1
$$
and
deduce
$$
\E(e^{\gamma L_{\ell_0}T})\le \kappa_5 (1+|x|^2+|y|^2).
$$
Since
$$
n_0=\inf\{ k,\;  V^{1,x,y}_{kT}=V^{2,x,y}_{kT}\}\le L_{\ell_0}+1
$$
it follows
$$
\E(e^{\gamma n_0T})\le  \kappa_5 (1+|x|^2+|y|^2)
$$
and
$$
\P(V^{1,x,y}_{kT}\ne V^{2,x,y}_{kT})=\P(k\ge n_0)\le \kappa_5(1+|x|^2+|y|^2) e^{-\gamma k T}.
$$
Moreover
$$
\begin{array}{ll}
\P(V^{1,x,y}_{kT+t}\ne V^{2,x,y}_{kT+t})\le \P(V^{1,x,y}_{kT}\ne V^{2,x,y}_{kT})&\le \kappa_5(1+|x|^2+|y|^2) e^{-\gamma k T}\\
& \le \kappa_6(1+|x|^2+|y|^2) e^{-\gamma (n T+t)}.
\end{array}
$$
We deduce for $\phi \in B_b(T^d)$ with $\sup_{x\in H}|\phi(x)|\le 1$
$$
\begin{array}{ll}
\left|\mathcal{P}_t[\phi](x) - \mathcal{P}_t[\phi](y)\right|&=\left|\E(\phi(V^{1,x,y}_{t})
-\phi(V^{2,x,y}_{t}))\right|\\
\\
&\le 2\kappa_6(1+|x|^2+|y|^2) e^{-\gamma t}.
\end{array}
$$

\qed

\subsection{Proof of the recurrence}\label{sec-app-rec}
Our method consists in applying  Proposition 3.4.5 in  \cite{DP2}, and
we apply Doob's Theorem (see Theorem 4.2.1 in  \cite{DP2}) in order to verify
the conditions of Proposition 3.4.5 in  \cite{DP2}.

\proof

Let us first introduce an auxiliary Markovian semigroup
 $\mathcal{R}_t[\phi](x)=\mathbb{E} \phi ( U^x_t )$
 corresponding to the Markov process $U$ where
 $$U^x_t=e^{tA} x + \int_0^te^{(t-s)A}G\,dW_s.$$
We denote
$$\mathcal{R}(t,x,\Gamma)= \mathbb{P}  ( U^x_t\in \Gamma)=
  \mathcal{R}_t[1_\Gamma](x),\qquad \Gamma \in \mathcal{B}(H),$$
the transition probabilities corresponding to $U$.

By Hypothesis $\ref{ipotesiuno}$, the Markovian semigroup
  $\mathcal{R}$ admits a unique invariant measure $\nu$ (see Theorem 6.3.3 in \cite{DP2}).
Moreover, $\mathcal{R}$ is irreducible (i.e., $\mathcal{R}(t,x,\Gamma)>0$ for all $\Gamma$ which is open and non empty set in $H$, see Theorem 7.2.1 in \cite{DP2}),
 and $t$-regular for any $t>0$ (i.e., the measures $\{\mathcal{R}(t,x,\cdot):
  x\in H\}$ are equivalent, see Theorem 7.3.1 in \cite{DP2}).

Next, recall that $\mathcal{P}_t[\phi](x)=\hat{\mathbb{E}} \phi (\hat X^x_t )$ (where $\hat{\mathbb{E}}$ means expectation with respect to probability $\hat{\mathbb{P}}$) is the
  Markovian semigroup  corresponding to weak solutions $\hat X^x$ of the equation (\ref{sde-coup})
  (we notice that the solutions of equation (\ref{sde-coup}) are unique in law). And we denote
 $$\mathcal{P}(t,x,\Gamma)= \hat{\mathbb{P}}  ( \hat{X}^x_t\in \Gamma )=
  \mathcal{P}_t[1_\Gamma](x),\qquad \Gamma \in \mathcal{B}(H)$$  the transition probabilities
  corresponding to $\hat{X}$. From the Girsanov theorem,
 the semigroup $\mathcal{P}_t$ can be represented, for $t\leq T$, by
 $$  \mathcal{P}_t[ \phi]
(x)=\mathbb{E}(\rho^{x}_T \phi(U^x_t)) $$ where $\rho^x_T=\exp\left(
-\int_0^t \<G^{-1}\Upsilon(U^x_t), dW_s\>_{\Xi}- \frac{1}{2}
\int_0^t |G^{-1}\Upsilon(U^x_t)|^2_{\Xi} ds\right) >0$,
$\mathbb{P}$-a.s.

Consequently, $\mathcal{P}$ is  irreducible (i.e., $\mathcal{P}(t,x,\Gamma)>0$ for all $\Gamma$ which is open and non empty set in $H$),
 and $t$-regular for any $t>0$ (i.e., the measures $\{\mathcal{P}(t,x,\cdot):
  x\in H\}$ are equivalent). And the above representation also implies that
 the semigroup $\mathcal{P}$ is stochastically continuous
 (i.e., $\lim_{t\searrow 0}  \mathcal{P}_t[\phi](x)=\phi(x)$ for all
   $x\in H$ and $\phi \in \mathcal{C}_b(H)$).

Moreover, by a similar argument to that of Theorem 8.4.4 in \cite{DP2},
there exists a
 measurable function $\eta: H\rightarrow \mathbb{R}^+$ such that
 $ \mu=\eta d\nu$ is the unique invariant measure corresponding to
 semigroup $(\mathcal{P}_t)_{t\ge 0}$. Indeed it is clear from the proof of
Theorem 8.4.3 in \cite{DP2} that Theorem 8.4.4 in \cite{DP2} remains true whenever $F$ can be approximated by a sequence of $\mathcal{C}^2_b$ functions converging in the bounded pointwise convergence sense.

$ $

We are now in  position to apply Doob's Theorem
(see Theorem 4.2.1 in \cite{DP2}) to obtain that

\begin{itemize}
 \item $\mu$ is strongly mixing and
 $\mathcal{P}_t(x,\Gamma)\rightarrow \mu(\Gamma)$ for all
  $\Gamma\in \mathcal{B}(H)$.
\item $\mu$ is equivalent to all measures $\mathcal{P}_t(x,\cdot)$.
\end{itemize}

Let $\Gamma$ be an open and non empty set in $H$. As  $\mu$ is equivalent to all measures $\mathcal{P}_t(x,\cdot)$ and $\mathcal{P}$ is  irreducible,
$\mu(\Gamma)>0$. Therefore,
$\liminf_{t\rightarrow\infty} \mathcal{P}_t(x,\Gamma)\rightarrow \mu(\Gamma)>0$. By
Proposition 3.4.5 in \cite{DP2}, the Markovian semigroup $\mathcal{P}$ is recurrent and
$\mathbb{P}\{\exists t>0 : \hat{X}_t\in \Gamma\}=1$.

In particular if, for all $T>0$, $\hat X^x$ is a weak solution of
equation (\ref{sde-coup}) in [0,T], under the probability
 $\hat{ \mathbb P} ^x_T$,
  and $\tau =\inf\{t : |\hat X^x_t|<
\epsilon\}$  then
$$\lim_{T\rightarrow \infty}\hat{ \mathbb P}^x_T\{ \tau^x<T \}=\lim_{T\rightarrow \infty}\hat{ \mathbb P}^x_T\{ \inf_{0\le t\le T}|\hat{X}_t^x|<\epsilon\}
=\lim_{T\rightarrow \infty}{ \mathbb P}\{ \inf_{0\le t\le T}|\hat{X}_t^x|<\epsilon\} =1. \qed$$

\end{document}